\documentclass[12pt,a4paper]{amsart}
\usepackage{amssymb,latexsym,graphicx}
\usepackage{color}
\usepackage[cp1252]{inputenc} 
\usepackage[T1]{fontenc}
\usepackage{lmodern}
\usepackage[english]{babel}
\usepackage[section]{placeins}
\usepackage[pagebackref=true]{hyperref}
\hypersetup{pdftex,colorlinks=true,linkcolor=blue,citecolor=blue,urlcolor=blue}
\usepackage{version}

\includeversion{weaker}

\includeversion{addendum}

\newtheorem{theorem}{Theorem}[section]
\newtheorem{lemma}[theorem]{Lemma}
\newtheorem{proposition}[theorem]{Proposition}

\newtheorem{claim}[theorem]{Claim}

\theoremstyle{definition} 

\theoremstyle{remark} 
\newtheorem{remark}[theorem]{Remark}

\numberwithin{equation}{section}
\numberwithin{figure}{section}
\numberwithin{table}{section}

\definecolor{purple}{RGB}{127,0,255}

\newenvironment{blue}{\color{blue}}{\ignorespacesafterend}

\newcommand{\N}{\mathbb{N}}

\newcommand{\R}{\mathbb{R}}

\newcommand{\Z}{\mathbb{Z}}

\newcommand{\noib}{\noindent $\bullet$~}
\newcommand{\noid}{\noindent $\diamond$~}
\newcommand{\ltop}[1]{\ell\downarrow #1}
\newcommand{\lbot}[1]{\ell\uparrow #1}
\newcommand{\page}[1]{\noindent \textbf{p.~#1}}

\newcommand{\da}{\frac{d}{dx}}
\newcommand{\db}[2]{\frac{d#1}{d#2}}
\newcommand{\ddb}[2]{\frac{d^2{#1}}{d{#2}^2}}

\DeclareMathOperator{\sign}{sign}
\newcommand{\pf}{\textbf{Proof.~}}
\newcommand{\bm}{\overline{m}}
\newcommand{\bN}{\overline{N}}

\begin{document}

\title[Sturm linear combinations eigenfunctions]{Sturm's theorem on zeros of linear combinations of eigenfunctions}

\author[P. B\'{e}rard]{Pierre B\'erard}
\author[B. Helffer]{Bernard Helffer}

\address{PB: Institut Fourier, Universit\'{e} Grenoble Alpes and CNRS, B.P.74\\ F38402 Saint Martin d'H\`{e}res Cedex, France.}
\email{pierrehberard@gmail.com}

\address{BH: Laboratoire Jean Leray, Universit\'{e} de Nantes and CNRS\\
F44322 Nantes Cedex, France.}
\email{Bernard.Helffer@univ-nantes.fr}


\date{October 3, 2017. Revised July 31, 2018}




\begin{abstract}
Motivated by recent questions about  the extension of Courant's nodal domain theorem, we revisit a theorem  published by C. Sturm in 1836, which  deals with zeros of linear combination of eigenfunctions of Sturm-Liouville problems. Although well known in the nineteenth century, this theorem seems to have been ignored or forgotten by some of the specialists in spectral theory since the second half of the twentieth-century. Although not specialists in History of Sciences, we have tried to put this theorem into the context of  nineteenth century mathematics.
\end{abstract}%

\maketitle

\begin{center}
\emph{To appear in Expositiones Mathematicae 2018, \\ except for the Appendices~\ref{S-weak} to \ref{S-cross} (in blue).\bigskip}
\end{center}%

\section{Introduction}\label{S-intro}

In this paper, we are interested in the following one-dimensional eigenvalue problem, where $r$ denotes the spectral parameter.

\begin{align}
& \da\left( K \db{V}{x}\right) + (r\, G - L) V = 0\,, \text{~for~} x \in ]\alpha , \beta[\,, \label{E-eq}\\[5pt]
& \left( K \db{V}{x} - h V\right)(\alpha)= 0\,, \label{E-bca} \\[5pt]
& \left( K \db{V}{x} + H V\right)(\beta)= 0\,. \label{E-bcb}
\end{align}

Here,
\begin{align}
& K, G, L : [\alpha,\beta] \to \R \text{~are positive functions}\,,\label{E-ass-0a}\\
& h\,, H \in [0,\infty] \text{~are non negative  constants, possibly infinite.}\label{E-ass-0b}
\end{align}

\begin{remark}\label{R-hH}
When $h=\infty$ (resp. $H=\infty$), the boundary condition should be understood as the Dirichlet boundary condition $V(\alpha)=0$ (resp. as the Dirichlet boundary condition $V(\beta)=0$).
\end{remark}%
\medskip

Precise assumptions on $K, G, L$ are given below.\medskip

Note that when $K=G\equiv 1$, \eqref{E-eq}--\eqref{E-bcb} is an eigenvalue problem for the classical operator $-\ddb{V}{x} + L V$.\medskip


This eigenvalue problem, in the above generality ($K, G, L$ functions of $x$), was first studied by Charles Sturm in a Memoir presented to the Paris Academy of sciences in September 1833, summarized in \cite{Sturm1833a,Sturm1833b}, and published in \cite{Sturm1836a,Sturm1836b}.

\begin{remark}\label{R-intro-0}
In this paper, we have mainly retained the notation of \cite{Sturm1836a}, except that we use $[\alpha,\beta]$ for the interval, instead of Sturm's notation $[\mathrm{x},\mathrm{X}]$. We otherwise use today notation and vocabulary. Note that in \cite{Sturm1836b}, Sturm uses lower case letters for the functions $K, G, L$, the same notation as Joseph Fourier in \cite{Fo1822}.
\end{remark}%


As far as the eigenvalue problem \eqref{E-eq}--\eqref{E-bcb} is concerned, Sturm's results can be roughly summarized in the following theorems.

\begin{theorem}[Sturm, 1836]\label{T-st1}
Under the assumptions \eqref{E-ass-0a}--\eqref{E-ass-0b}, the ei\-gen\-value problem \eqref{E-eq}--\eqref{E-bcb} admits an increasing infinite sequence $\{\rho_i, i \ge 1\}$ of positive simple eigenvalues, tending to infinity. Furthermore, the associated eigenfunctions $V_i$ have the following remarkable property: the function $V_i$ vanishes, and changes sign, precisely $(i-1)$ times in the open interval $]\alpha,\beta[\,$.
\end{theorem}%

\begin{theorem}[Sturm, 1836]\label{T-st2r}
Let $Y = A_m V_m + \cdots + A_n V_n$ be a  non trivial linear combination of eigenfunctions of the eigenvalue problem \eqref{E-eq}--\eqref{E-bcb}, with $1 \le m \le n$, and $\{A_j, m \le j \le n\}$ real constants such that $A_m^2 + \cdots + A_n^2 \not = 0$. Then, the function $Y$ has at least $(m-1)$, and at most $(n-1)$ zeros in the open interval $]\alpha,\beta[$.
\end{theorem}%

The first theorem today appears in most textbooks on Sturm-Liouville theory. Although well known in the nineteenth century, the second theorem (as well as the more precise Theorem~\ref{T-st2}) seems to have been ignored or forgotten by some of the specialists in spectral theory since the second half of the twentieth-century, as the following chronology indicates.\medskip

\begin{description}
  \item[1833] Sturm's Memoir presented to the Paris Academy of sciences in September, summarized in \cite{Sturm1833a,Sturm1833b}.
  \item[1836] Sturm's papers \cite{Sturm1836a,Sturm1836b} published. Joseph Liouville summarizes Sturm's results in \cite[$\S$~{III}, p.~257]{Liou1836a}, and uses them to study the expansion of a given function $f$ into a series of eigenfunctions of \eqref{E-eq}.
  \item[1877]  Lord Rayleigh writes \emph{``a beautiful theorem has been discovered by Sturm''} as he mentions Theorem~\ref{T-st2r} in \cite[Section~142]{Ray1877}. 
  \item[1891] F. Pockels \cite[pp. 68-73]{Po1891} gives a summary of Sturm's results, including Theorem~\ref{T-st2r}, and mentions the different proofs provided by Sturm, Liouville and Rayleigh. On the basis of a note of Sturm in F\'{e}russac's Bulletin \cite{Sturm1829}, Pockels (p.~71, lines 12-17) also suggests that Sturm may have looked for a statement in higher dimension as well, without success. Sturm indeed mentions studying an example with spherical symmetry in dimension $3$ (leading to an ordinary differential equation with singularity), to which he may have applied Theorem~\ref{T-st2r}.
  \item[1903] Hurwitz \cite{Hu1903} gives a lower bound for the number of zeros of the sum of a trigonometric series with a spectral gap and refers, somewhat inaccurately, to Sturm's Theorems. This result, known as the Sturm-Hurwitz theorem, already appears in a more general framework in Liouville's paper \cite{Liou1836a}. \\
      See \cite[$\S$~2]{ErNo2004} for a generalization of the Sturm-Hurwitz theorem to Fourier integrals with a spectral gap, \cite{OvTa2005} for geometric applications, and the recent paper \cite{Ste2018} which quantifies the Sturm-Hurwitz theorem.
  \item[1931] Courant \& Hilbert \cite{CH1931,CH1953} extensively  mention the Sturm-Liouville problem. They do not refer to the original papers of  Sturm, but to B\^{o}cher's book \cite{Bo1917}  which does not include Theorem~\ref{T-st2r}. They then state an extension of the so-called Courant's nodal domain theorem to linear combination of eigenfunctions, \cite[footnote, p.~394]{CH1931} and \cite[footnote, p.~454]{CH1953}, and refer to the dissertation of H.~Herrmann \cite{H1932}. It turns out that neither Herrmann's dissertation, nor his later papers, consider this extension of Courant's Theorem.
 \item[1950] The book \cite{GaKr2002} by F.~Gantmacher and M.~Krein contains several notes on Sturm's contributions. One result (Corollary, Chap.~{III.5}, p.~138),  stated in the context of Chebyshev systems, is stronger than Theorem~\ref{T-st2r}, yet weaker than Theorem~\ref{T-st2}. The book does however not mention \cite{Sturm1836b}.
  \item[1956]  Pleijel mentions Sturm's Theorem~\ref{T-st2r}, somewhat inaccurately, in \cite[p.~543 and 550]{Pl}.
  \item[1973] V.~Arnold \cite{Arn1973} points out that an extension of Courant's theorem to linear combinations of eigenfunctions cannot be true in general. Counter\-examples were first given  by O.~Viro for the $3$-sphere (with the canonical metric)  \cite{Vir1979} and, more recently in the papers \cite{BH2017,BH2018a}, see also \cite{GZ2003}.\\
It seems to us that Arnold may have not been aware of Theorem~\ref{T-st2r}. Indeed, in \cite{Arn2011},  see also the Supplementary problem 9 in \cite[p.~327]{Arn1992}, he mentions a proof, suggested by I.~Gelfand, of the upper bound in Theorem~\ref{T-st2r}. Gelfand's idea is to  ``use fermions rather than bosons'', and  to apply Courant's nodal domain theorem in the fermionic context. However, Arnold concludes by writing \cite[p.~30]{Arn2011}, ``\emph{the arguments {\rm[given by Gelfand]} do not yet provide a proof}''. It is interesting to note that  Liouville's and Rayleigh's proofs of the lower bound in Theorem~\ref{T-st2r} use an idea similar to Gelfand's, see the proof of Claim~\ref{C-lip-2}.\\
As far as we know, the first implementation of Gelfand's idea into a complete proof of Theorem~\ref{T-st2r} is given in \cite{BH2018a,BH2018b}.
\end{description}

\begin{remark}\label{R-intro-GaHa}
In \cite{Sturm1836b}, Theorem~\ref{T-st2r}  first appears as a corollary to a much deeper theorem  \cite[$\S$~{XXIV}]{Sturm1836b}, in which Sturm describes the time evolution of the $x$-zeros of a solution $u(x,t)$ of the heat equation. We shall not consider this topic here, and we refer to \cite{GaHa2005,LuMi2009} for modern formulations and a historical analysis.
\end{remark}%

Our interest in Theorem~\ref{T-st2r} arose from reading \cite{Ku}, and investigating Courant's nodal domain theorem and its extension to linear combination of eigenfunctions.\medskip

The main purpose of this paper is to popularize Theorems~\ref{T-st2r} and \ref{T-st2}, as well as Sturm's originality and ideas. Sturm's results are clearly stated in the summaries \cite{Sturm1833a,Sturm1833b}. Unfortunately, Sturm's detailed papers \cite{Sturm1836a,Sturm1836b} are written linearly, and contain very few tagged statements. Our second purpose is to provide an accessible proof of Theorems~\ref{T-st2r} and \ref{T-st2}, meeting  today standards of rigor. We in particular state precise assumptions, clarify some technical points, and provide some alternative proofs. We otherwise closely follow the original proofs, and we provide precise cross-references to Sturm's papers.
\medskip

In this paper, we make the following strong assumptions.
\begin{equation}\label{E-intro-ass-S}
\left\{
\begin{array}{l}
[\alpha,\beta] \subset ]\alpha_0,\beta_0[\,,\\[5pt]
K, G, L \in C^{\infty}(]\alpha_0,\beta_0[)\,,\\[5pt]
K, G, L >0 \text{~on~} ]\alpha_0,\beta_0[\,.
\end{array}%
\right.
\end{equation}

\begin{remark}\label{R-intro-12}
Neither Sturm nor Liouville make any explicit regularity assumptions, see Subsection~\ref{SS-hist-ass} and Remark~\ref{R-weaker} for more details.
\end{remark}%

\subsection*{Organization of the paper}

In Section~\ref{S-stp}, we prove Theorem~\ref{T-st2}, Sturm's refined version of Theorem~\ref{T-st2r}, following the ideas of \cite[$\S$~XXVI]{Sturm1836b}. In Section~\ref{S-lip}, we prove Theorem~\ref{T-LiSt}, Liouville's version of Theorem~\ref{T-st2r}, following \cite{Liou1836a,Liou1836b}. In Section~\ref{S-hist}, we describe the context of Sturm's papers and his ideas. Appendix~\ref{S-lim} provides the detailed proof of a technical argument. Appendix~\ref{S-Trans} contains the citations from Sturm's papers in their original French formulation. Appendix~\ref{S-weak} considers Sturm's theorem under weaker assumptions. Appendices~\ref{S-sop} and \ref{S-cross} provide cross-references between our paper and the papers of Sturm and Liouville.

\subsection*{Acknowledgements} The authors would like to thank N.~Kuznetsov and J.~L\"{u}tzen for their comments on an earlier version of this paper. They also thank the anonymous referees for their constructive remarks.


\section{Sturm's o.d.e. proof of Theorem~\ref{T-st2r}}\label{S-stp}

\subsection{Preliminary lemmas and notation}\label{SS-stp-pre}

\subsubsection{} Recall that $\{(\rho_j,V_j), j\ge 1\}$ are the eigenvalues and eigenfunctions of the eigenvalue problem \eqref{E-eq}--\eqref{E-bcb}.\medskip

By our assumption $L >0$, the eigenvalues are positive, $\rho_j > 0$.
Under the Assumptions~\eqref{E-intro-ass-S}, the functions $V_j$ are $C^{\infty}$ on $]\alpha_0,\beta_0[\,$.   This follows from Cauchy's existence and uniqueness theorem, or from Liouville's existence proof \cite{Liou1836a}. Note that the assumption $L > 0$ is convenient, but not necessary. It actually suffices that $\frac{L}{G}$ be bounded from below.\medskip

In this section, we fix
\begin{equation}\label{E-stp-2a}
Y = \sum_{j=m}^{n} A_j V_j\,,
\end{equation}
a linear combination of eigenfunctions of the eigenvalue problem \eqref{E-eq}--\eqref{E-bcb}, where $1 \le m \le n\,$, and where the $A_j$ are real constants.

\begin{remark}\label{R-st2r-1}
We shall always assume that $Y \not \equiv 0\,$, which is equivalent to assuming that $\sum_m^n A_j^2 \neq 0\,$. As far as the statement of Theorem~\ref{T-st2r} is concerned, and without loss of generality, it is simpler to assume that $A_m \, A_n \neq 0\,$.
\end{remark}%


We also introduce the associated family of functions, $\{Y_k, k\in \Z\}$, where
\begin{equation}\label{E-stp-2}
Y_k = (-1)^k\, \sum_{j=m}^n \rho_j^k\, A_j \, V_j\,.
\end{equation}
Note that $Y_0$ is the original linear combination~$Y$, and that $Y_k \equiv 0$ if and only if $Y \equiv 0\,$. \medskip

Roughly speaking, Sturm's idea is to show that the number of zeros of $Y_k$, in the interval $]\alpha,\beta[$, is non-decreasing with respect to $k$, and then to take the limit when $k$ tends to infinity, see Subsection~\ref{SS-stp}. Up to changing the constants $A_j$, it suffices to compare the numbers of zeros of $Y$ and $Y_1$. For this purpose, Sturm compares the signs of $Y$ and $Y_1$ near the zeros of $Y$ (Lemma~\ref{L-stp-3}), and at the non-zero local extrema of $Y$ (Lemmas~\ref{L-stp-6} and \ref{L-stp-8}). The main ingredient for this purpose is the differential relation \eqref{E-stp-6k}.  In the sequel, we indicate the pages in Sturm's papers corresponding to the different steps of the proof.

\subsubsection{} For $m \le p \le n$, write the equations satisfied by the eigenfunction $V_p$,
\begin{align}\label{E-stp-4p}
& \da\left( K \db{V_p}{x}\right) + (\rho_p\, G - L) V_p = 0\,,\\
& \left( K \db{V_p}{x} - h V_p\right)(\alpha)= 0\,, \\
& \left( K \db{V_p}{x} + H V_p\right)(\beta)= 0\,,
\end{align}
and multiply the $p$-th equation by $\rho_p^k\, A_p$.  Summing up from $p=m$ to $n$, yields the following lemma.

\begin{lemma}\label{L-stp-2}
Assume that \eqref{E-intro-ass-S} holds. Let $k \in \Z$.
\begin{enumerate}
  \item The function $Y_k$ satisfies the boundary conditions \eqref{E-bca} and \eqref{E-bcb}.
  \item The functions $Y_k$ and $Y_{k+1}$ satisfy the  differential relation
\begin{equation}\label{E-stp-6k}
G\, Y_{k+1} = K\, \ddb{Y_k}{x} + \db{K}{x}\, \db{Y_k}{x} - L\, Y_k\,.
\end{equation}
  \item  Under the Assumptions~\eqref{E-intro-ass-S}, the function $Y_k$ cannot vanish at infinite order at a point $\xi \in [\alpha,\beta]$,  unless $Y \equiv 0\,$.
\end{enumerate}
\end{lemma}%

\pf  \cite[p.~437]{Sturm1836b}~ Assertions (1) and (2) are clear by linearity.\medskip

For Assertion (3), assume that  $Y_k \not \equiv 0\,$, and that it  vanishes at infinite order at some $\xi$. Then, according to \eqref{E-stp-6k} and its successive derivatives, the function $Y_{k+1}$ also vanishes at infinite order at $\xi$, and so does $Y_{\ell}$ for any $\ell \ge k$.  Assume, as indicated in Remark~\ref{R-st2r-1}, that $A_n \not = 0$.  Fixing some $p \ge 0$, we can write, for any $\ell \ge k$,
\begin{equation*}
\frac{d^pV_n}{dx^p}(\xi) + \sum_{j=m}^{n-1} \left( \frac{\rho_j}{\rho_n}\right)^{\ell} \frac{A_j}{A_n} \, \frac{d^pV_j}{dx^p}(\xi) = 0\,.
\end{equation*}

Since $\rho_n > \rho_j$ for $m \le j \le n-1$, letting $\ell$ tend to infinity, we conclude that $\frac{d^pV_n}{dx^p}(\xi) = 0$. This would be true for all $p$, which is impossible by Cauchy's uniqueness theorem, or by Sturm's argument \cite[$\S$~{II}]{Sturm1836a}. \hfill \qed

\begin{remark}\label{R-stp-L2}
Assertion (3), and the fact that the zeros of $Y$ are isolated, with finite multiplicities, are implicit in \cite{Sturm1836b}.
\end{remark}%

\begin{lemma}\label{L-stp-3}
Assume that \eqref{E-intro-ass-S} holds. Let $U$ denote any $Y_k$, and $U_1 = Y_{k+1}$. Let $\xi \in [\alpha,\beta]$ be a zero of $U$, of order $p \ge 2\,$. Then, there exist constants $B_{\xi}$ and $B_{1,\xi}\,$, and smooth functions $R_{\xi}$ and $R_{1,\xi}\,$, such that
\begin{equation}\label{E-stp-L3}
\left\{
\begin{array}{l}
  U(x) =  B_{\xi} (x-\xi)^p + (x-\xi)^{p+1}R_{\xi}(x)\,,\\[5pt]
  U_{1}(x) =  B_{1,\xi} (x-\xi)^{p-2} + (x-\xi)^{p-1}R_{1,\xi}(x)\,,\\[5pt]
\text{with~}  B_{\xi} \, B_{1,\xi}  > 0\,.
\end{array}
\right.
\end{equation}
\end{lemma}%

\pf \cite[p.~439]{Sturm1836b}~ Assume that $\xi$ is a zero of order $p \ge 2$ of $U$, so that
$$
U(\xi) = \cdots = \frac{d^{p-1}U}{dx^{p-1}}(\xi)= 0
$$
and
$$
\frac{d^pU}{dx^p}(\xi) \not = 0\,.
$$

Taylor's formula with integral remainder term,  see Laplace \cite[p.~179]{Lapl1820} (in Livre premier, Partie 2, Chap. 3, \S~44), gives the existence of some function $R_{\xi}$ such that
\begin{equation*}
U(x) = B_{\xi} (x-\xi)^p + (x-\xi)^{p+1}R_{\xi}(x)\,,
\end{equation*}
where $$ B_{\xi} = \frac{1}{p!}\frac{d^pU}{dx^p}(\xi) \not = 0\,.$$
Equation \eqref{E-stp-6k} implies that
\begin{equation*}
(G U_{1})(x) = p(p-1) B_{\xi} (x-\xi)^{p-2}K(x) + (x-\xi)^{p-1}S_{\xi}(x)\,,
\end{equation*}
for some smooth function $S_{\xi}$. It follows that
$$
U_{1}(x) = B_{1,\xi} (x-\xi)^{p-2} + (x-\xi)^{p-1}R_{1,\xi}(x)\,,
$$
for some function $R_{1,\xi}$\,, with $B_{1,\xi} = p(p-1) \frac{K(\xi)}{G(\xi)}B_{\xi}\,$.\\
 In particular, $B_{1,\xi} \, B_{\xi} > 0\,$ and this proves the lemma. \hfill \qed\medskip

\begin{lemma}\label{L-stp-4}
Assume that \eqref{E-intro-ass-S} holds. Assume that $h \in [0,\infty[\,$, i.e., that the boundary condition at $\alpha$ is not the Di\-ri\-chlet boundary condition.
Let $U$ denote any $Y_k$, $U_1 = Y_{k+1}$, and assume that $U(\alpha) = 0\,$.  Then, $\alpha$ is a zero of $U$ of even order, i.e., there exists $n_U \in \N\setminus \{0\}$ such that $\frac{d^pU}{dx^p}(\alpha)=0$ for $0 \le p \le 2n_U-1$ and $\not = 0$ for $p = 2n_U$.\\[5pt]
When $H\in [0,\infty[\,$, a similar statement holds at the boundary $\beta$.
\end{lemma}%

\pf \cite[p.~440-441]{Sturm1836b}~ Assume that $U(\alpha) = 0$. By Lemma~\ref{L-stp-2}, $U$ does not vanish at infinite order at $\alpha$, so that there exists $p\geq 1$ with $$U(\alpha) = \cdots = \frac{d^{p-1}Y_k}{dx^{p-1}}(\alpha)= 0$$ and $$\frac{d^pU}{dx^p}(\alpha) \not = 0\,.$$ Taylor's formula with integral remainder term gives
\begin{equation*}
U(x) = B_{\alpha} (x-\alpha)^p + (x-\alpha)^{p+1}R_{\alpha}(x)\,,
\end{equation*}
where $B_{\alpha} = \frac{1}{p!}\frac{d^pU}{dx^p}(\alpha) \not = 0\,$. \\
The boundary condition at $\alpha$ implies that $\db{U}{x}(\alpha) = 0$, and hence that $p \ge 2$. By Lemma~\ref{L-stp-3}, we can write
$$U_{1}(x) = B_{1,\alpha} (x-\alpha)^{p-2} + (x-\alpha)^{p-1}R_{1,\alpha}(x)\,,$$ with $B_{1,\alpha} \, B_{\alpha} > 0\,$.\\
If $p=2$, then $U_{1}(\alpha) \not = 0$. If $p > 2$, one can continue.\\
If $p = 2q$, one arrives at $$Y_{k+q}(x) = B_{k+q,\alpha} + (x-\alpha)R_{k+q,\alpha}(x)\,,$$
 with $Y_{k+q}(\alpha) = B_{k+q,\alpha}$ and $B_{k+q,\alpha} \, B_{k,\alpha} > 0\,$.\\
If $p = 2q+1$, one arrives at $$ Y_{k+q}(x) = B_{k+q,\alpha}(x-\alpha) + (x-\alpha)^2 R_{k+q,\alpha}(x)\,,$$
with $B_{k+q\alpha} \, B_{k,\alpha} > 0$ and $\db{Y_{k+q}}{x}(\alpha) = B_{k+q,\alpha} \not = 0\,$. On the other-hand, since $Y_{k+q}$ satisfies \eqref{E-bca} and $Y_{k+q}(\alpha)=0\,$, we must have $\db{Y_{k+q}}{x}(\alpha) = 0\,$,  because $h > 0\,$. This yields a contradiction and proves that the case $p = 2q+1$ cannot occur. The lemma is proved. \hfill \qed\medskip

\subsection{Counting zeros}\label{SS-stp-cz}

Assume that \eqref{E-intro-ass-S} holds. Let $U$ denote any $Y_k$, and $U_1 = Y_{k+1}$. They satisfy the relation \eqref{E-stp-6k}.\medskip

From Lemma~\ref{L-stp-2}, we know that $U$ cannot vanish at infinite order at a point $\xi \in [\alpha,\beta]$. If $\xi \in ]\alpha,\beta[$ and $U(\xi) = 0\,$, we define the \emph{multiplicity} $m(U,\xi)$ of the zero $\xi$ by
\begin{equation}\label{E-stp-10}
m(U,\xi) = \min\{p ~|~\frac{d^pU}{dx^p}(\xi) \not = 0 \}\,.
\end{equation}

From Lemma~\ref{L-stp-4}, we know that the multiplicity $m(U,\alpha)$ is even if $h \in [0,\infty[$, and that the multiplicity $m(U,\beta)$ is even if $H \in [0,\infty[\,$. We define the \emph{reduced multiplicity} of $\alpha$ by
\begin{equation}\label{E-stp-12}
\bm (U,\alpha) = \left\{
\begin{array}{ll}
\frac{1}{2}m(U,\alpha) &\text{if~} h \in [0,\infty[\,,\\[5pt]
0 &\text{if~} h = \infty\,,
\end{array}%
\right.
\end{equation}
and a similar formula for the reduced multiplicity of $\beta$.\medskip

By Lemma~\ref{L-stp-2}, the function $U$ has finitely many distinct zeros \break $\xi_1(U) < \xi_2(U) < \cdots < \xi_p(U)$ in the interval $]\alpha,\beta[\,$. We define the \emph{number of zeros of} $U$ \emph{in} $]\alpha,\beta[\,$, \emph{counted with multiplicities}, by
\begin{equation}\label{E-stp-14}
N_m(U,]\alpha,\beta[) = \sum_{j=1}^p m(U,\xi_i(U))\,,
\end{equation}
and we use the notation $N_m(U)$ whenever the interval is clear.\medskip

We define the \emph{number of zeros of} $U$ \emph{in} $[\alpha,\beta]$, \emph{counted with multiplicities}, by
\begin{equation}\label{E-stp-16}
\bN _m(U,[\alpha,\beta]) = \sum_{j=1}^p m(U,\xi_i(U)) + \bm (U,\alpha) + \bm (U,\beta)\,,
\end{equation}
and we use the notation $\bN _m(U)$ whenever the interval is clear.\medskip

We define the \emph{number of zeros of} $U$ \emph{in} $]\alpha,\beta[$ (\emph{multiplicities not accounted for}) by
\begin{equation}\label{E-stp-17}
N(U,]\alpha,\beta[) = p\,,
\end{equation}
and we use the notation $N(U)$ whenever the interval is clear.\medskip

Finally, we define the \emph{number of sign changes of} $U$ \emph{in the interval} $]\alpha,\beta[$ by
\begin{equation}\label{E-stp-18}
N_v(U,]\alpha,\beta[) = \sum_{j=1}^p \frac{1}{2}\left[ 1 - (-1)^{m(U,\xi_j(U))}\right]\,.
\end{equation}

\begin{remark}\label{R-cs}
Note that sign changes of the function $Y$ correspond to zeros with odd multiplicity.
\end{remark}

\subsection{Comparing the numbers of zeros of $Y_k$ and $Y_{k+1}$}\label{SS-stp}

Assume that \eqref{E-intro-ass-S} holds. Let $U$ be some $Y_k$ and $U_1 = Y_{k+1}$. In this subsection, we show that the number of zeros of $U_1$ is not smaller than the number of zeros of $U$.\medskip

\begin{lemma}\label{L-stp-6}
Let $\xi < \eta$ be two zeros of $U$ in $[\alpha,\beta]$. Then, there exists some $a_{\xi,\eta} \in ]\xi,\eta[$ such that $U(a_{\xi,\eta}) \, U_1(a_{\xi,\eta}) < 0\,$.
\end{lemma}%

\begin{remark}\label{R-Lstpf-6}
We do not assume that $\xi, \eta$ are consecutive zeros.
\end{remark}%

\pf \cite[p.~437]{Sturm1836b}~ Since $U$ cannot vanish identically in $]\xi,\eta[$ (see Lemma~\ref{L-stp-2}), there exists some $x_0 \in ]\xi,\eta[$ such that $U(x_0) \neq 0$. Let $\varepsilon_0 = \sign (U(x_0))$. Then $\varepsilon_0 U$ takes a positive value at $x_0$, and hence $M := \sup\{\varepsilon_0 U(x) ~|~ x \in [\xi,\eta]\}$ is positive and achieved at some $a_{\xi,\eta} \in ]\xi,\eta[$. Denote this point by $a$ for short, then, $$\varepsilon_0 U(a) > 0\,, \,\db{U}{x}(a) = 0 \mbox{  and } \varepsilon_0 \ddb{U}{x}(a) \le 0\,.
$$
 It follows from \eqref{E-stp-6k} that $\varepsilon_0 U_1(a) < 0\,$, or equivalently, that \break $U(a) U_1(a) < 0\,$. The lemma is proved. \hfill \qed

\begin{lemma}\label{L-stp-8}
Let $\xi \in ]\alpha,\beta]$. Assume that $U(\xi) = 0\,$, and that $U$ does not change sign in $]\alpha,\xi[$. Then, there exists some $a_{\xi} \in [\alpha,\xi[$ such that $U(a_{\xi}) U_1(a_{\xi}) < 0$.\\
Let $\eta \in [\alpha,\beta[$. Assume that $U(\eta) = 0$, and that $U$ does not change sign in $]\eta,\beta[$. Then, there exists some $b_{\eta} \in ]\eta,\beta]$ such that $U(b_{\eta}) U_1(b_{\eta}) < 0$.
\end{lemma}%

\pf \cite[p.~438]{Sturm1836b}~ Since $U$ cannot vanish identically in $]\alpha,\xi[$ (see Lemma~\ref{L-stp-2}), there exists $x_0 \in ]\alpha,\xi[$ such that $U(x_0) \neq 0\,$. Let $\varepsilon_{\xi} = \sign (U(x_0))\,$. Since $U$ does not change sign in $]\alpha,\xi[\,$, $\varepsilon_{\xi} U(x) \ge 0$ in $]\alpha,\xi[\,$. Then,
$$
M_{\xi} := \sup \{\varepsilon_{\xi}U(x) ~|~ x \in [\alpha,\xi]\} > 0\,.
$$
Let
$$
a_{\xi} := \inf \{x \in [\alpha,\xi] ~|~ \varepsilon_{\xi} U(x) = M_{\xi}\}\,.
$$
Then $a_{\xi} \in [\alpha,\xi[$.\medskip

If $a_{\xi} \in ]\alpha,\xi[$, then $\varepsilon_{\xi} U(a_{\xi}) > 0\,$, $\db{U}{x}(a_{\xi}) = 0\,$, and $\varepsilon_{\xi} \ddb{U}{x}(a_{\xi}) \le 0\,$. By \eqref{E-stp-6k}, this implies that $\varepsilon_{\xi} U_1(a_{\xi}) < 0$. Equivalently, $ U(a_{\xi}) U_1(a_{\xi}) < 0$.\medskip

\begin{claim}\label{C-Lstp-8}
If $a_{\xi} = \alpha\,$, then $\varepsilon_{\xi} \,U(\alpha) > 0\,$, $h=0\,$, $\db{U}{x}(\alpha) = 0\,$, and $\varepsilon_{\xi} \, \ddb{U}{x}(\alpha) \le 0\,$.
\end{claim}%

\emph{Proof of the claim.~} Assume that $a_{\xi} = \alpha\,$, then $\varepsilon_{\xi} \, U(\alpha) > 0\,$, and hence $h \neq \infty\,$. If $h$  were in $]0,\infty[$, we would have $\varepsilon_{\xi}\, \db{U}{x}(\alpha) = h \varepsilon_{\xi} U(\alpha) > 0\,$, and hence $a_{\xi} > \alpha\,$. It follows that the assumption $a_{\xi} = \alpha$ implies that $h=0$ and $\db{U}{x}(\alpha) = 0$. If $\varepsilon_{\xi} \ddb{U}{x}(\alpha)$ where positive, we would have $a_{\xi} > \alpha$. Therefore, the assumption $a_{\xi} = \alpha$ also implies that $\varepsilon_{\xi} \ddb{U}{x}(\alpha) \le 0\,$. The claim is proved.\medskip

If $a_{\xi} = \alpha$, then by Claim~\ref{C-Lstp-8} and \eqref{E-stp-6k}, we have $\varepsilon_{\xi} U_1(\alpha) < 0$. Equivalently, $U(\alpha) U_1(\alpha) < 0$. The first assertion of the lemma is proved. The proof of the second assertion is similar. \hfill \qed


\begin{proposition}\label{P-stp-10}
Assume that \eqref{E-intro-ass-S} holds, and let $k \in \Z$. Then,
\begin{equation}\label{E-stp-24}
N_v(Y_{k+1},]\alpha,\beta[) \ge N_v(Y_k,]\alpha,\beta[)\,,
\end{equation}
i.e., in the interval $]\alpha,\beta[$, the function $Y_{k+1}$ changes sign at least as many times as the function $Y_k$.
\end{proposition}%

\pf \cite[p.~437-439]{Sturm1836b}~~ We keep the notation $U = Y_k$ and $U_1 = Y_{k+1}$.  By Lemma~\ref{L-stp-2}, the functions $U$ and $U_1$ have finitely many zeros in $]\alpha,\beta[$, with finite multiplicities. Since $\alpha$ and $\beta$ are fixed, we skip the mention to the interval $]\alpha,\beta[$ in the proof, and we examine several cases.\medskip

\emph{Case~1}. If $N_v(U) = 0$, there is nothing to prove.\medskip

\emph{Case~2}. Assume that $N_v(U)=1$. Then $U$ admits a unique zero $\xi \in ]\alpha,\beta[$ having odd multiplicity. Without loss of generality, we may assume that $U \ge 0$ in $]\alpha,\xi[$ and $U \le 0$ in $]\xi,\beta[\,$.  By Lemma~\ref{L-stp-8}, there exist $a \in [\alpha,\xi[$ and $b \in ]\xi,\beta]$ such that $U_1(a) < 0$ and $U_1(b) > 0\,$. \medskip

It follows that the function $U_1$ vanishes and changes sign at least once in $]\alpha,\beta[\,$, so that $N_v(U_1) \ge 1 = N_v(U)$, which proves the lemma in Case~2.
\medskip

\emph{Case~3}. If $N_v(U) = 2$, the function $U$ has exactly two zeros, having odd multiplicities, $\xi$ and $\eta$ in $]\alpha,\beta[\,$, $\alpha < \xi < \eta < \beta\,$, and we may assume that $U|_{]\alpha,\xi[} \ge 0\,$,
$U|_{]\xi,\eta[} \le 0$, and $U|_{]\eta,\beta[} \ge 0$. The arguments given in Case~2 imply that there exist $a \in [\alpha,\xi[$ such that $U_1(a) < 0$ and $b \in ]\eta,\beta]$ such that $U_1(b) < 0\,$. In $]\xi,\eta[$ the function $U$ does not vanish identically and therefore achieves a global minimum at a point $c$ such that $U(c) < 0\,$, $\db{U}{x}(c)=0\,$, and $\ddb{U}{x}(c) \ge 0\,$. Equation\eqref{E-stp-6k} then implies that $U_1(c) >0\,$. \medskip

We can conclude that the function $U_1$ vanishes and changes sign at least twice in $]\alpha,\beta[$, so that $N_v(U_1) \ge 2 = N_v(U)$.
\medskip

\emph{Case~4}. Assume that $N_v(U) = p \ge 3\,$. Then, $U$ has exactly $p$ zeros, with odd multiplicities, in $]\alpha,\beta[\,$, $\alpha < \xi_1 < \xi_2 < \cdots < \xi_p < \beta$\,, and one can assume that
\begin{equation*}
\begin{array}{l}
U|_{]\alpha,\xi_1[} \ge 0\,, \,(-1)^p U|_{]\xi_p,\beta[} \ge 0\,, \mbox{~and~}\\[5pt]
(-1)^iU|_{]\xi_i,\xi_{i+1}[} \ge 0 \mbox{ for } 1 \le i \le p-1\,.
\end{array}%
\end{equation*}%
One can repeat the arguments given in the Cases~2 and 3, and conclude that there exist $a_0, \ldots, a_p$ with $a_0 \in [\alpha,\xi_1[\,$, $a_i \in ]\xi_i,\xi_{i+1}[$ for $1 \le i \le p-1$, and $a_p \in ]\xi_p,\beta]$ such that $(-1)^iU_1(a_i) < 0\,$. \medskip

We can then conclude that the function $U_1$ vanishes and changes sign at least $p$ times in $]\alpha,\beta[\,$, i.e. that $N_v(U_1) \ge p = N_v(U)\,$. \medskip

This concludes the proof of Proposition~\ref{P-stp-10}.\hfill \qed \medskip

\begin{proposition}\label{P-stp-12}
Assume that \eqref{E-intro-ass-S} holds. For any $k \in \Z$,
\begin{equation}\label{E-stp-26}
N_m(Y_{k+1},]\alpha,\beta[) \ge N_m(Y_{k},]\alpha,\beta[)\,,
\end{equation}
 i.e., in the interval $]\alpha,\beta[\,$, counting multiplicities of zeros, the function $Y_{k+1}$ vanishes at least as many times as the function $Y_{k}$.
\end{proposition}%

\pf  \cite[p.~439-442]{Sturm1836b}~ Let $U=Y_k$ and $U_1=Y_{k+1}\,$.  If $U$ does not vanish in $]\alpha,\beta[\,$, there is nothing to prove. We now assume that $U$ has at least one zero in $]\alpha,\beta[\,$. By Lemma~\ref{L-stp-2}, $U$ and $U_1$ have finitely many zeros in $]\alpha,\beta[\,$. Let
$$
\alpha < \xi_1 < \cdots < \xi_k < \beta
$$
be the distinct zeros of $U$, with multiplicities $p_i = m(U,\xi_i)$ for $1 \le i \le k\,$. Let $\sigma_0$ be the sign of $U$ in $]\alpha,\xi_1[\,$, $\sigma_i$ the sign of $U$ in $]\xi_i,\xi_{i+1}[$ for $1 \le i \le k-1\,$, and $\sigma_k$ the sign of $U$ in $]\xi_k,\beta[\,$. Note that
$$
\sigma_i = \sign\left( \frac{d^{p_i}U}{dx^{p_i}}(\xi_i)\right) \text{~for~} 1\le i \le k\,.
$$

By Lemma~\ref{L-stp-8}, there exist $a_0 \in [\alpha,\xi_1[$ and $a_k \in ]\xi_k,\beta]$ such that $U(a_0) U_1(a_0) < 0$ and $U(a_k) U_1(a_k) < 0\,$.
By Lemma~\ref{L-stp-6}, there exists $a_i \in ]\xi_i,\xi_{i+1}[\,$, $1 \le i \le k-1\,$, such that $U(a_i) U_1(a_i) < 0\,$. \\
Summarizing, we have obtained:
\begin{equation}\label{E-stp-P122}
\text{For~} 0 \le i \le k\,, ~~ U(a_i) \,U_1(a_i) < 0 \,.
\end{equation}

We have the relation
\begin{equation}\label{E-stp-P124}
N_m(U,]\alpha,\beta[) = \sum_{i=1}^k N_m(U,]a_{i-1},a_i[) = \sum_{i=1}^k p_i\,.
\end{equation}
Indeed, for $1 \le i \le k$, the interval $]a_{i-1},a_i[$ contains precisely one zero $\xi_i$ of $U$, with multiplicity $p_i\,$.

For $U_1$, we have the inequality
\begin{equation}\label{E-stp-P126}
N_m(U_1,]\alpha,\beta[) \ge \sum_{i=1}^k N_m(U_1,]a_{i-1},a_i[)\,,
\end{equation}
because $U_1$ might have zeros in the interval $]\alpha,a_0[$ if $a_0 > \alpha$ (resp. in the interval $]a_k,\beta[$ if $a_k < \beta$).

\begin{claim}\label{C-stp-P12}
For $1 \le i \le k$, $$N_m(U_1,]a_{i-1},a_i[) \ge N_m(U,]a_{i-1},a_i[) = p_i\,.$$
\end{claim}%

To prove the claim, we consider several cases.\medskip

\noib If $p_i = 1$\,, then $U(a_{i-1}) \, U(a_i) < 0$ and, by \eqref{E-stp-P122}, $U_1(a_{i-1}) U_1(a_i) < 0\,$, so that $N_m(U_1,]a_{i-1},a_i[) \ge 1\,$.\medskip

\noib If $p_i \ge 2$, we apply Lemma~\eqref{L-stp-3} at $\xi_i\,$: there exist real numbers $B, B_1$ and smooth functions $R$ and $R_{1}$, such that, in a neighborhood of $\xi_i\,$,
\begin{equation}\label{E-stp-P128}
\left\{
\begin{array}{ll}
U(x) &= B(x-\xi_i)^{p_i} + (x-\xi_i)^{p_i+1}R(x)\,,\\[5pt]
U_1(x) &= B_1(x-\xi_i)^{p_i-2} + (x-\xi_i)^{p_i-1}R_{1}(x)\,,
\end{array}%
\right.
\end{equation}
where $\sign(B) = \sign(B_1) = \sigma_i\,$.

We now use \eqref{E-stp-P122} and the fact that $\sign(U(a_i)) = \sigma_i\,$.\medskip

\noid If $p_i\ge 2$ is odd, then $\sigma_{i-1} \sigma_i = -1\,$. It follows that
$$
\sigma_i U_1(a_i) < 0 \text{~and~} \sigma_i U_1(a_{i-1}) > 0\,.
$$
By \eqref{E-stp-P128}, for $\varepsilon$ small enough, we also have
$$
\sigma_i U_1(\xi_i+\varepsilon) > 0 \text{~and~} \sigma_i U_1(\xi_i-\varepsilon) < 0 \,.
$$
This means that $U_1$ vanishes at order $p_i-2$ at $\xi_i\,$, and at least once in the intervals $]a_{i-1},\xi_i-\varepsilon[$ and $]\xi_i+\varepsilon,a_i[$, so that
$$
N_m(U_1,]a_{i-1},a_i[) \ge p_i -2 + 2 = p_i = N_m(U,]a_{i-1},a_i[)\,.
$$

\noid If $p_i\ge 2$ is even, then $\sigma_{i-1} \sigma_i = 1\,$. It follows that
$$
\sigma_i U_1(a_i) < 0 \text{~and~} \sigma_i U_1(a_{i-1}) < 0\,.
$$
By \eqref{E-stp-P128}, for $\varepsilon$ small enough, we also have
$$
\sigma_i U_1(\xi_i+\varepsilon) > 0 \text{~and~} \sigma_i U_1(\xi_i-\varepsilon) > 0 \,.
$$
This means that $U_1$ vanishes at order $p_i-2$ at $\xi_i\,$, and at least once in the intervals $]a_{i-1},\xi_i-\varepsilon[$ and $]\xi_i+\varepsilon,a_i[\,$, so that
$$
N_m(U_1,]a_{i-1},a_i[) \ge p_i -2 + 2 = p_i = N_m(U,]a_{i-1},a_i[)\,.
$$

The claim is proved, and the proposition as well. \hfill \qed

\begin{proposition}\label{P-stp-14}
Assume that \eqref{E-intro-ass-S} holds. For any $k \in \Z$,
\begin{equation}\label{E-stp-28}
\bN _m(Y_{k+1},[\alpha,\beta]) \ge \bN _m(Y_{k},[\alpha,\beta])\,,
\end{equation}
i.e., in the interval $[\alpha,\beta]$, counting multiplicities of interior zeros, and reduced multiplicities of $\alpha$ and $\beta$, the function $Y_{k+1}$ vanishes at least as many times as the function $Y_{k}$.
\end{proposition}%

\pf \cite[p.~440-442]{Sturm1836b}~ Recall that the reduced multiplicity of $\alpha$ (resp. $\beta$) is zero if the Dirichlet condition holds at $\alpha$ (resp. at $\beta$) or if $U(\alpha) \neq 0$ (resp. $U(\beta) \neq 0$). Furthermore, according to Lemma~\ref{L-stp-4}, if $h \in [0,\infty[$ and $U(\alpha)=0$ (resp. if $H\in [0,\infty[$ and $U(\beta) = 0$), then $m(U,\alpha) = 2p$ (resp. $m(U,\beta)=2q$). \medskip

\emph{Case~1.~} Assume that $N_m(U,]\alpha,\beta[) = 0\,$. Without loss of generality, we may assume that $U > 0$ in $]\alpha,\beta[\,$.\medskip

\noib If $U(\alpha) \neq 0$ and $U(\beta) \neq 0$, there is nothing to prove.\medskip

\noib Assume that $U(\alpha) = U(\beta) = 0\,$. Then, there exists $a \in ]\alpha,\beta[$ such that $$
 U(a) = \sup\{U(x) ~|~ x \in[\alpha,\beta]\}\,,$$
  with
  $$ U(a) > 0\,, \db{U}{x}(a)=0\,, \mbox{ and }\ddb{U}{x}(a) \le 0\,.
 $$
 It follows from \eqref{E-stp-6k} that $U_1(a) < 0\,$, and that
\begin{equation}\label{E-stp-P144}
\bN_m(U_1,[\alpha,\beta] )= \bN(U_1,[\alpha,a]) + \bN(U_1,[a,\beta])\,.
\end{equation}
It now suffices to look separately at the intervals $[\alpha,a]$ and
$[a,\beta]\,$.

\noid Interval $[\alpha,a]$. If the Dirichlet condition holds at $\alpha$, there is nothing to prove. If $h \in [0,\infty[$, $m(U,\alpha) = 2p \ge 2$ and, by Lemma~\ref{L-stp-3},
\begin{equation}\label{E-stp-P146}
\begin{array}{l}
U(x) = B (x-\alpha)^{2p} + (x-\alpha)^{2p+1}R(x)\,,\\[5pt]
U_1(x) = B_1 (x-\alpha)^{2p-2} + (x-\alpha)^{2p-1}R_1(x)\,,\\[5pt]
\text{with~} B > 0 \text{~and~} B_1 > 0\,.
\end{array}%
\end{equation}
It follows that $U_1(\alpha+\varepsilon) > 0$ for any positive $\varepsilon$ small enough so that $N_m(U_1,]\alpha,a[) \ge 1$. It follows that
\begin{equation}\label{E-stp-P148}
\begin{array}{ll}
\bN_m(U_1,[\alpha,a]) & = \bm(U_1,\alpha) + N_m(U_1,]\alpha,a[)\ge p-1+1\,,\\[5pt]
\text{i.e.}\\[5pt]
\bN_m(U_1,[\alpha,a]) & \ge \bN(U,[\alpha,a])\,.
\end{array}%
\end{equation}

\noid Interval $[a,\beta]$. The proof is similar.\medskip

\noib Assume that $U(\alpha) = 0$ and $U(\beta) \neq 0\,$. The proof is similar to the previous one with $a \in ]a,\beta]\,$.\medskip

\noib Assume that $U(\alpha) \neq 0$ and $U(\beta) = 0\,$. The proof is similar to the previous one with $a \in [\alpha,a[\,$. \medskip

\emph{Case~2.~} Assume that $N_m(U,]\alpha,\beta[) \ge 1\,$. \medskip

\noib If $U(\alpha) \neq 0$ (resp. $U(\beta) \neq 0$), there is nothing to prove for the boundary $\alpha$ (resp. $\beta$).  \medskip

\noib If $U(\alpha)=0$ (resp. $U(\beta) =0$), the number $a_0$ (resp. $a_k$) which appears in the proof of Proposition~\ref{P-stp-12} belongs to the open interval $]\alpha,\xi_1[$ (resp. to the open interval $]\xi_k,\beta[$), where $\xi_1$ (resp. $\xi_k$) is the smallest (resp. largest) zero of $U$ in $]\alpha,\beta[$. We can then apply the proof of Step.~1 to the interval $[\alpha,a_0]$ (resp. to the interval $[a_k,\beta]$) to prove that $\bN_m(U_1,[\alpha,a_0]) \ge \bN_m(U,[\alpha,a_0]$) (resp. to prove that $\bN_m(U_1,[a_k,\beta]) \ge \bN_m(U,[a_k,\beta]$). This proves Proposition~\ref{P-stp-14}. \hfill \qed \medskip

We can now state Sturm's refined version of Theorem~\ref{T-st2r}.


\begin{theorem}\label{T-st2}
Assume that \eqref{E-intro-ass-S} holds, and let $Y$ be the non trivial linear combination
\begin{equation}\label{E-st2-a}
Y = \sum_{p=m}^n A_p V_p\,,
\end{equation}
where $1 \le m \le n$, and where $\{A_p, m\le p \le n\}$ are real constants such that $A_m^2 + \cdots + A_n^2 \not \equiv 0\,$.
Then, with the notation of Subsection~\ref{SS-stp-cz},
\begin{equation}\label{E-st2-b}
N_v(Y,]\alpha,\beta[)  \le N_m(Y,]\alpha,\beta[) \le \bN_m(Y,[\alpha,\beta])\,,
\end{equation}
\begin{equation}\label{E-st2-s}
(m-1)  \le N_v(Y,]\alpha,\beta[)\,
 \text{~and~} \, \bN_m(Y,[\alpha,\beta])  \le (n-1)\,.
\end{equation}
\end{theorem}%

\pf \cite[p.~442]{Sturm1836b}~ Let $N(V)$ be any of the above functions. We may of course assume that $A_m \not = 0$ and $A_n \not = 0$. In the preceding lemmas, we have proved that $N(Y_{k+1}) \ge N(Y_k)$ for any $k \in \Z$. This inequality can also be rewritten as
\begin{equation}\label{E-st2-ba}
N(Y_{(-k)}) \le N(Y) \le N(Y_k) \text{~~for any~~} k \ge 1\,.
\end{equation}

Letting $k$ tend to infinity, we conclude that
\begin{equation}\label{E-st2-c}
N(V_m) \le N(Y) \le N(V_n)\,,
\end{equation}
and we can apply Theorem~\ref{T-st1}.\hfill \qed \medskip

\textbf{Remark}. For a complete proof of the limiting argument when $k$ tends to infinity, we refer to Appendix~\ref{S-lim}.

\section{Liouville's approach to Theorem~\ref{T-st2r}}\label{S-lip}

\subsection{Main statement}
We keep the notation of Section~\ref{S-stp}. Starting from a linear combination $Y$ as in \eqref{E-stp-2a},  Liouville also considers the family $Y_k$ given by \eqref{E-stp-2}, and shows that the number of zeros of $Y_{k+1}$ is not smaller than the number of zeros of $Y_k$. His proof is based on a generalization of Rolle's theorem.

\begin{remark}\label{R-lip-0}
In his proof, Liouville \cite{Liou1836b} only considers the zeros in the open interval $]\alpha,\beta[\,$.
\end{remark}%

As in Section~\ref{S-stp}, for $1 \le m \le n$, we fix $Y = \sum_{j=m}^n A_j V_j$, a linear combination of eigenfunctions of the eigenvalue problem \eqref{E-eq}--\eqref{E-bcb}, and we assume that $A_m A_n \neq 0\,$, see Remark~\ref{R-st2r-1}.

\begin{theorem}\label{T-LiSt}
Counting zeros with multiplicities in the interval $]\alpha,\beta[\,$, the function $Y$ (1) has at most $(n-1)$ zeros
and, (2) has at least $(m-1)$ zeros.
\end{theorem}%

\pf Liouville uses the following version of Rolle's theorem  (Michel Rolle (1652-1719) was a French mathema\-ti\-cian). This version of Rolle's theorem seems to go back to Cauchy and Lagrange.

\begin{lemma}\label{L-lip-2}
Let $f$ be a  function in $]\alpha_0,\beta_0[\,$. Assume that $$f(x') = f(x'') = 0 \mbox{ for some  }x', x''\,, \,\alpha_0 < x' < x'' < \beta_0\,.$$
\begin{enumerate}
 \item If the  function $f$ is differentiable, and has $\nu - 1$ distinct zeros in the interval $]x',x''[\,$, then the  derivative $f'$ has at least $\nu$ distinct zeros in $]x',x''[\,$.
  \item If the  function $f$ is smooth, and has $\mu-1$ zeros counted with multiplicities in the interval $]x',x''[\,$, then the derivative $f'$ has at least $\mu$ zeros counted with multiplicities in $]x',x''[\,$.
 \end{enumerate}
\end{lemma}%

\textbf{Proof of the lemma.~} Call $x_1 < x_2 < \cdots x_{\nu-1}$ the distinct zeros of $f$ in $]x',x''[\,$. Since $f(x')=f(x'')=0\,$,  by Rolle's theorem \cite{Rolle1691}, the function $f'$ vanishes at least once in each open interval determined by the $x_j\,$, $1 \le j \le \nu-1\,$, as well as in the intervals $]x',x_1[$ and $]x_{\nu-1},x''[\,$. It follows that $f'$ has at least $\nu$ distinct zeros in $]x',x''[\,$, which proves the first assertion.\medskip

Call $m_j$ the multiplicity of the zero $x_j$, $1 \le j \le \nu - 1\,$. Then $f'$ has at least $\nu$ zeros, one in each of the open intervals determined by $x', x''$ and the $x_j$'s, and has a zero at each $x_j$ with multiplicity $m_j-1\,$, provided that $m_j > 1$. It follows that the number of zeros of $f'$ in $]x',x''[\,$, counting multiplicities, is at least
$$
\sum_{j=1}^{\nu - 1}(m_j - 1) + \nu = \sum_{j=1}^{\nu - 1}m_j + 1\,,
$$
which proves the second assertion. \hfill \qed \medskip

\subsection{Proof of the assertion~~\emph{``$Y$ has at most $(n-1)$ zeros in $]\alpha,\beta[$, counting multiplicities''}}\label{SS-lip-1}~\\[5pt]
Write \eqref{E-eq} for $V_1$ and for $V_p$, for some $m \le p \le n$. Multiply the first equation by $-V_p$, the second by $V_1$, and add the resulting equations. Then
\begin{equation}\label{E-lip-2}
V_1 \da\left( K \db{V_p}{x}\right) - V_p \da\left( K \db{V_1}{x}\right) + (\rho_p - \rho_1) G V_1 V_p = 0 \,.
\end{equation}

Use the identity
\begin{equation}\label{E-lip-4}
V_1 \da\left( K \db{V_p}{x}\right) - V_p \da\left( K \db{V_1}{x}\right) = \da \left( V_1 K \db{V_p}{x} - V_p K \db{V_1}{x}\right)\,,
\end{equation}
and integrate from $\alpha$ to $t$ to get the identity
\begin{equation}\label{E-lip-6}
(\rho_1 - \rho_p) \int_{\alpha}^t G V_1 V_p \, dx = K(t) \left( V_1(t) \db{V_p}{x}(t) - V_p(t) \db{V_1}{x}(t)\right)\,.
\end{equation}
Here we have used  the boundary condition \eqref{E-bca} which implies that $$\left( V_1(\alpha) \db{V_p}{x}(\alpha) - V_p(\alpha) \db{V_1}{x}(\alpha)\right)=0\,.$$
Multiplying the identity \eqref{E-lip-6} by $A_p$, and summing for $p$ from $m$ to $n$, we obtain
\begin{equation}\label{E-lip-8}
\int_{\alpha}^t  G V_1  \sum_{p=m}^n (\rho_1 - \rho_p) A_p V_p \, dx = K(t) \left( V_1 \db{Y}{x} - Y \db{V_1}{x}\right)(t)\,.
\end{equation}
or
\begin{equation}\label{E-lip-10}
\int_{\alpha}^t  G V_1  \sum_{p=m}^n (\rho_1 - \rho_p) A_p V_p \, dx = K(t)\, V_1^2(t) \,\frac{d}{dt}\left( \frac{Y}{V_1}\right)(t)\,,
\end{equation}
 where we have used the fact that the function $V_1$ does not vanish in the interval $]\alpha,\beta[\,$.\medskip

Let $\Psi (x) = \frac{Y}{V_1}(x)$. The zeros of $Y$ in $]\alpha,\beta[$ are the same as the zeros of $\Psi$, with the same multiplicities. Let $\mu$ be the number of zeros of $Y$, counted with multiplicities. Using  Lemma~\ref{L-lip-2}, Assertion~(2), one can show that $\db{\Psi}{x}$ has at least $\mu - 1$ zeros in $]\alpha,\beta[\,$, and hence so does the left-hand side of \eqref{E-lip-10},
\begin{equation*}
\int_{\alpha}^t  G V_1  \sum_{p=m}^n (\rho_1 - \rho_p) A_p V_p \, dx\,.
\end{equation*}
On the other hand, this function vanishes at $\alpha$ and $\beta$ (because of the boundary condition \eqref{E-bcb} or orthogonality). By Lemma~\ref{L-lip-2}, its derivative,
\begin{equation}\label{E-lip-12}
 V_1  \sum_{p=m}^n (\rho_1 - \rho_p) A_p V_p
\end{equation}
has at least $\mu$ zeros counted with multiplicities in $]\alpha,\beta[\,$. We have proved the following

\begin{lemma}\label{L-lip-4}
If the function $Y = \sum_{p=m}^n A_p V_p$ has  at least $\mu$ zeros counted with multiplicities in the interval $]\alpha,\beta[\,$, then the function \break
$Y_1 = \sum_{p=m}^n (\rho_1 - \rho_p) A_p V_p$ has at least $\mu$ zeros, counted with multiplicities, in $]\alpha,\beta[\,$.
\end{lemma}%

Applying this lemma iteratively, we deduce that if $Y$ has  at least $\mu$ zeros counted with multiplicities in $]\alpha,\beta[\,$, then, for any $k \ge 1$, the function
\begin{equation}\label{E-lip-14}
Y_k = \sum_{p=m}^n (\rho_1 - \rho_p)^k A_p V_p
\end{equation}
has at least $\mu$ zeros, counted with multiplicities,  in $]\alpha,\beta[\,$.\medskip

We may of course assume that the coefficient $A_n$ is non-zero. The above assertion can be rewritten as the statement:\medskip

{\it For all $k \ge 0\,$, the equation
\begin{equation}\label{E-lip-16}
A_m \left( \frac{\rho_m - \rho_1}{\rho_n - \rho_1} \right)^k V_m + \cdots + A_{n-1} \left( \frac{\rho_{n-1} - \rho_1}{\rho_n - \rho_1} \right)^k V_{n-1} + A_n V_n = 0
\end{equation}
has at least $\mu$ solutions in $]\alpha,\beta[$, counting multiplicities. }\medskip

Letting $k$ tend to infinity, and using the fact that $V_n$ has exactly $(n-1)$ zeros in $]\alpha,\beta[\,$, this implies that $\mu \le (n-1)$. This proves the first assertion. \hfill \qed

\subsection{Proof of the assertion~~ \emph{``$Y$ has at least $(m-1)$ zeros in $]\alpha,\beta[$, counting multiplicities''}}\label{SS-lip-2}~\\[5pt]
We have seen that the number of zeros of $Y_k$ is less than or equal to the number of zeros of the function $Y_{k+1}$. This assertion actually holds for any $k \in \Z$, and can also be rewritten as,
\begin{equation}\label{E-lip-30}
N_m(Y_{-k}) \le N_m(Y)\,,
\end{equation}
for any $k \ge 0$, where
\begin{equation}\label{E-lip-30a}
Y_{-k} = A_m (\rho_m - \rho_1)^{-k} V_m + \cdots + A_n (\rho_n - \rho_1)^{-k} V_n\,,
\end{equation}
and we can again let $k$ tend to infinity. The second assertion is proved and Theorem \ref{T-LiSt} as well.\hfill \qed

\subsection{Liouville's 2nd approach to the 2nd part of Theorem~\ref{T-LiSt}}\label{SS-lip-3}
If the function $Y$ has $\mu_1$ distinct zeros, and $\mu \le \mu_1$ sign changes, we call $a_i$, $\alpha < a_1 < \cdots < a_{\mu} <\beta$, the points at which $Y$ changes sign.

\begin{claim}\label{C-lip-2}
The function $Y$ changes sign at least $(m-1)$ times in the interval $]\alpha,\beta[\,$.
\end{claim}%

\textbf{Proof of the claim.} Assume,  by contradiction,  that $\mu \le (m-2)$. Consider the function
\begin{equation}\label{E-lip-32}
x \mapsto W(x) := \Delta(a_1,\ldots,a_{\mu};x)\,,
\end{equation}
where the function $\Delta$ is defined as the determinant
\begin{equation}\label{E-lip-34}
\begin{vmatrix}
  V_1(a_1) & V_1(a_2) & \cdots & V_1(a_{\mu}) & V_1(x) \\
  V_2(a_1) & V_2(a_2) & \cdots & V_2(a_{\mu}) & V_2(x)\\
  \vdots & \vdots & \vdots & \vdots & \vdots\\
  V_{\mu+1}(a_1) & V_{\mu+1}(a_2) & \cdots & V_{\mu+1}(a_{\mu}) & V_{\mu+1}(x)\\
\end{vmatrix}
\,.
\end{equation}

The function $W$ vanishes at the points $a_i\,, 1 \le i \le \mu\,$. According to the first part in Theorem~\ref{T-LiSt}, $W$ being a linear combination of the first $\mu+1$ eigenfunctions, vanishes at most $\mu$ times in $]\alpha,\beta[\,$, counting multiplicities. This implies that each zero $a_i$ of $W$ has order one, and that $W$ does not have any other zero in $]\alpha,\beta[\,$. It follows that the function $Y W$ vanishes only at the points $\{a_i\}$, $1 \le i \le \mu$, and that it does not change sign. We can assume that $YW \ge 0\,$. On the other hand, we have
\begin{equation}\label{E-lip-36}
\int_{\alpha}^{\beta} G Y W \, dx = 0\,,
\end{equation}
because $Y$ involves the functions $V_p$ with $p \ge m$ and $W$ the functions $V_q$ with $q \le \mu + 1 \le m-1\,$.  This gives a contradiction. \hfill \qed

\begin{remark}\label{R-lip-E140}
Liouville does actually not use the determinant \eqref{E-lip-34}, but a similar approach, see \cite[p.~259]{Liou1836a}, Lemme~$1^{\text{er}}$.  The determinant $\Delta$ appears in \cite[Section~142]{Ray1877}. The paper \cite{BH2018b} is based on a careful analysis of this determinant.
\end{remark}%

\begin{remark}\label{R-weaker}
The arguments in Subsection~\ref{SS-lip-1}, using Assertion~(1) of
Lemma~\ref{L-lip-2}, instead of Assertion~(2), yield an upper bound
on the number  of zeros of $Y$, multiplicities not accounted for. This estimate holds under weaker regularity assumptions, namely only assuming that the functions $G, L$ are continuous, and that the function $K$ is $C^1$, see Appendix~\ref{S-weak}, and compare with \cite{GaKr2002}, Chap.~{III.5}.
\end{remark}%

\section{Mathematical context of Sturm's papers.\\
Sturm's motivations and ideas}\label{S-hist}

\subsection{On Sturm's style}

Sturm's papers \cite{Sturm1836a,Sturm1836b} are written in French, and  quite long, about 80 pages each. One difficulty in reading them is the lack of layout structure. The papers are written linearly, and divided into sequences of sections, without any title. Most results are stated without tags, ``Theorem'' and the like, and only appear in the body of the text. For example, \cite{Sturm1836a} only contains one theorem stated as such, see $\S$~{XII}, p.~125. In order to have an overview of the results contained in \cite{Sturm1836a}, the reader should look at the announcement \cite{Sturm1833a}. Theorem~\ref{T-st2r} is stated in \cite{Sturm1833b}.\medskip

For a more thorough analysis of Sturm's papers on differential equations, we refer to \cite{LuMi2009,GaHa2005}.  We refer to \cite{Bo1911, Sina1999} for the relationships between Theorem~\ref{T-st1} and Sturm's theorem on the number of real roots of real polynomials.

\subsection{Sturm's motivations}

Sturm's motivations come from mathematical physics, and more precisely, from the problem of heat diffusion in a non-homogeneous bar. He considers the heat equation,
\begin{equation}\label{E-hist-2}
G \frac{\partial u}{\partial t} = \frac{\partial}{\partial x}\left( K \frac{\partial u}{\partial x} \right) -L u\,, \text{~for~} (x,t) \in ]\alpha,\beta[\times \R_{+}\,,
\end{equation}
with boundary conditions
\begin{equation}\label{E-hist-2bc}
\left\{
\begin{array}{l}
K(\alpha)\, \frac{\partial u}{\partial x}(\alpha,t) - h\, u(\alpha,t) = 0\,,\\[5pt]
K(\beta)\, \frac{\partial u}{\partial x}(\beta,t) + H\, u(\beta,t) = 0\,,
\end{array}%
\right.
\end{equation}
for all $t > 0$, and with the initial condition
\begin{equation}\label{E-hist-2ic}
u(x,0) = f(x)\,, \text{~for~} x \in ]\alpha,\beta[\,,
\end{equation}
where $f$ is a given function.\medskip

The functions $K, G, L$ and the constants $h, H$ describe the physical properties of the bar, see \cite[Introduction, p.~376]{Sturm1836b}. Sturm refers to the book of Sim\'{e}on Denis Poisson \cite{Pois1835}, rather than to Fourier's book \cite{Fo1822}, because Poisson's equations are more general, see \cite[Chap.~{III}]{Sina1999}. \medskip

The boundary conditions \eqref{E-bca}-\eqref{E-bcb} and \eqref{E-hist-2bc} first appeared in the work of Fourier \cite{Fo1822} but are called ``Robin's condition'' in the recent literature. Victor Gustave Robin (1855-1897) was a French mathematician.\medskip

As was popularized by Fourier  and Poisson,  in order to solve \eqref{E-hist-2}, Sturm uses the method of separation of variables, and is therefore led to the eigenvalue problem \eqref{E-eq}--\eqref{E-bcb}.\medskip

\subsection{ Sturm's assumptions}\label{SS-hist-ass}

In \cite{Sturm1836a,Sturm1836b}, Sturm  implicitly as\-su\-mes that the functions $K, G, L$ are $C^{\infty}$ and, explicitly, that $K$ is positive, see \cite[p.~108]{Sturm1836a}. For the eigenvalue problem, he also assumes that $G, L$ are positive, see \cite[p.~381]{Sturm1836b}. In \cite[p.~394]{Sturm1836b}, he mentions that $L$ could take negative values, and implicitly assumes, in this case, that $\frac{L}{G}$ is bounded from below. \medskip

In \cite{Liou1836b}, Liouville does not mention any regularity assumption on the functions $G, K, L$.  He however indicates a regularity assumption (piecewise $C^2$ functions) in a previous paper,
\cite[Footnote~$(*)$, p.~256]{Liou1836a}.

\subsection{Sturm's originality}

Before explaining Sturm's proofs, we would like to insist on the \emph{originality} of his approach. Indeed, unlike his predecessors, Sturm does not look for explicit solutions of the differential  equation \eqref{E-eq-1}  (i.e., solutions in closed form, or given as sums of series or as integrals), but he rather looks for \emph{qualitative properties} of the solutions, properties which can be deduced directly from the differential equation itself. The following excerpts are translated from \cite[Introduction]{Sturm1836a}\footnote{See Appendix~\ref{S-Trans} for the original citations in French.}.\\
{\it
One only knows how to integrate these equations in a very small number of particular cases, and one can otherwise not even obtain a first integral; even when one knows the expression of the function which satisfies such an equation, in finite form, as a series, as integrals either definite or indefinite, it is most generally difficult to recognize in this expression the behaviour and the characteristic properties of this function. \ldots \\
Although it is important to be able to determine the value of the unknown function for an isolated value of the variable it depends upon, it is not less necessary to discuss the behaviour of this function, or otherwise stated, the form and the twists and turns of the curve whose ordinate would be the function, and the abscissa the independent variable. It turns out that one can achieve this goal by the sole consideration of the differential equation themselves, without having to integrate them. This is the purpose of the present memoir. ~\ldots
}

\subsection{Sturm and the existence and uniqueness theorem for ordinary differential equation}

In \cite[p.~108]{Sturm1836a}, Sturm considers the differential equation
\begin{equation*}\label{E-sturm-I}
\da \left( K \db{V}{x} \right) + GV = 0\,, \hspace{4cm} (I)
\end{equation*}
and takes the existence and uniqueness theorem for gran\-ted. More precisely, he claims \cite[p.~108]{Sturm1836a}, without any reference whatsoever,\\
\emph{The complete integral of equation (I) must contain two arbitrary constants, for which one can take the values of $V$ and of $\db{V}{x}$ corresponding to some particular value of $x$. Once these values are fixed, the function $V$ is fully determined by equation (I), it has a uniquely determined value for each value of $x$.}\\
On the other hand, he gives two  arguments for the fact that a solution of (I) and its derivative cannot vanish simultaneously at a point without vanishing identically, see \cite[$\S$~{II}]{Sturm1836a}.  When the coefficients $K, G$ of the differential equation depend upon a parameter $m$, e.g. continuously, Sturm also takes for granted the fact that the solution $V(x,m)$, and its zeros, depend continuously on $m$.\medskip

In \cite[$\S$~{II}]{Sturm1836b}, Sturm mentions the existence proof given by Liouville in \cite{Liou1836a}, see also \cite{Liou1830}. According to \cite{Gila1989}, Augustin-Louis Cauchy may have presented the existence and uniqueness theorem for ordinary differential equations in his course at \'{E}cole polytechnique as early as in the year 1817-1818. Following a recommendation of the administration of the school, Cauchy delivered the notes of his lectures in 1824, see \cite{Cauc1824} and, in particular, the introduction by Christian Gilain who discovered these notes in  1974. These notes apparently had a limited distribution. Liouville entered the \'{E}cole polytechnique in 1825, and there attended the mathematics course given by Amp\`{e}re\footnote{We are grateful to J.~L\"{u}tzen for providing this information.} (as a matter of fact Amp\`{e}re and Cauchy gave the course every other year, alternatively). Liouville's proof of the existence theorem for differential equations in \cite{Liou1830}, \`{a} la Picard but before Picard, though limited to the particular case of 2nd order linear equations, might be the first well circulated proof  of an existence theorem for differential equations, see \cite[$\S$~34]{Lutz1984}. Cauchy's theorem was later popularized in the second volume of Moigno's book, published in 1844, see \cite{Moig1840}, ``Vingt-sixi\`{e}me Le\c{c}on'' $\S$~159, pp. 385--396.

\subsection{Sturm's proof of Theorem~\ref{T-st1}}

Theorem~\ref{T-st1} is proved in \cite{Sturm1836b}. For the first assertion, see $\S$~{III} (p.~384) to {VII}; for the second assertion, see $\S$~{VIII} (p.~396) to {X}. \medskip

The proof is based on the paper \cite{Sturm1836a} in which Sturm studies the zeros of the solution of the initial value problem,
\begin{align}
& \da\left( K(x,m) \db{V}{x}(x,m)\right) + G(x,m) V(x,m) = 0\,,\label{E-eq-1}\\
& \left( K \db{V}{x} - h V\right)(\alpha,m)= 0\,. \label{E-bca-1}
\end{align}
Here $K,G$ are assumed to be functions of $x$ depending on a real parameter $m$, with $K$ positive (the constants $h$ and $H$ may also depend on the parameter $m$). The solution $V(x,m)$ is well defined up to a scaling factor. The main part of \cite{Sturm1836a} is devoted to studying how the zeros of the function $V(x,m)$ (and other related functions) depend on the parameter $m$, see \cite[$\S$~{XII}, p.~125]{Sturm1836a}. While developing this program, Sturm proves the \emph{oscillation, separation and comparison theorems} which nowadays bear his name, \cite[$\S$~{XV}, {XVI} and {XXXVII}]{Sturm1836a}. \medskip

The eigenvalue problem \eqref{E-eq}--\eqref{E-bcb} itself is studied in \cite{Sturm1836b}. For this purpose, Sturm considers the functions
$$
K(x,r) \equiv K(x) \text{~and~ }G(x,r) = r G(x) - L(x)\,,
$$
the solution $V(x,r)$ of the corresponding initial value problem \eqref{E-eq-1}--\eqref{E-bca-1}, and applies the results and methods of \cite{Sturm1836a}.\medskip

The spectral data of the eigenvalue problem \eqref{E-eq}--\eqref{E-bcb} are determined by the following transcendental equation in the spectral parameter $r$,
\begin{equation}\label{E-te}
K(\beta) \db{V}{x}(\beta,r) + H V(\beta,r) = 0\,,
\end{equation}
see, \cite{Sturm1836b}, \S{III}, page 383, line 8 from bottom.

\subsection{Sturm's two proofs of Theorem~\ref{T-st2r}}

Theorem~\ref{T-st2r} appears in \cite[$\S$~{XXV}, p.~431]{Sturm1836b}, see also the announcement \cite{Sturm1833b}. \medskip

Sturm's general motivation, see the introductions to \cite{Sturm1836a} and \cite{Sturm1836b}, was the investigation of heat diffusion in a (non-homogeneous) bar, whose physical properties are described by the functions $K, G, L$. He first obtained Theorem~\ref{T-st2r}  as a corollary of a much deeper theorem which describes the behaviour, as time varies, of the $x$-zeros of a solution $u(x,t)$ of the heat equation \eqref{E-hist-2}-\eqref{E-hist-2ic}. When the initial temperature $u(x,0)$ is given by a linear combination of simple states,
\begin{equation}\label{E-stp-n2}
u(x,0) = Y(x) = \sum_{j=m}^n A_j V_j
\end{equation}
the function $u(x,t)$ is given by
\begin{equation}\label{E-stp-n4}
u(x,t) = \sum_{j=m}^{n} e^{-t\rho_j}A_j V_j\,.
\end{equation}
When $t$ tends to infinity, the $x$-zeros of $u(x,t)$ approach those of $V_p$, where $p$ is the least integer $j, m\le j \le n$ such that $A_j \not = 0$.\medskip

J.~Liouville, who was aware of  Theorem~\ref{T-st2r}, made  use of it in \cite{Liou1836a}, and provided a purely  ``ordinary differential equation'' proof in \cite{Liou1836b}, a few months before the actual publication of \cite{Sturm1836b}. This induced Sturm to provide two proofs of Theorem~\ref{T-st2r} in \cite{Sturm1836b}, his initial proof using the heat equation, and another proof based on the sole ordinary differential equation. The proofs of Sturm actually give a more precise result. In \cite[p.~379]{Sturm1836b}, Sturm writes,\\
{\it
M. Liouville gave a direct proof of this theorem, which for me was a mere corollary of the preceding one, without taking care of the particular case in which the function vanishes at one of the extremities of the bar. I have also found, after him, another direct proof which I give in this memoir. M.~Liouville made use of the same theorem in a very nice memoir which he published in the July issue of his journal, and which deals with the expansion of an arbitrary function into a series made of the functions $V$ which we have considered.
}

\medskip

The time independent analog to studying the behaviour of the $x$-zeros of \eqref{E-stp-n4} is to study the behaviour of the zeros of the family of functions $\{Y_k\}_{k \in \Z}$, where
\begin{equation}\label{E-stp-nk}
Y_k(x) = \sum_{j=m} \rho_j^k \, A_j  V_j\,,
\end{equation}
as $k$ tends to infinity.\bigskip

\appendix

\section{The limiting argument in~\eqref{E-lip-16}}\label{S-lim}

Recall that we assume that $A_n \neq 0$. Define
\begin{equation}\label{E-lip-20}
\omega = \left( \frac{\rho_{n-1} - \rho_1}{\rho_n - \rho_1} \right)^k \,.
\end{equation}

One can rewrite \eqref{E-lip-16} as
$$
V_n(x) + \omega \, \Pi(x) = 0 \,,
$$ where
\begin{equation}\label{E-lip-22}
\Pi(x) = \sum_{p=m}^{n-1} \left( \frac{\rho_p - \rho_1}{\rho_{n-1} - \rho_1} \right)^k \frac{A_p}{A_n} \, V_p \,.
\end{equation}

It follows that $\Pi$ is uniformly bounded by
\begin{equation}\label{E-lip-24}
\left| \Pi(x) \right| \le M := n \, \max_p \left| \frac{A_p}{A_n} \right| \, \max_p \sup_{[\alpha,\beta]}|V_p|\,.
\end{equation}
Similarly,
\begin{equation}\label{E-lip-26}
\left| \db{\Pi}{x}(x) \right| \le N := n \, \max_p \left| \frac{A_p}{A_n} \right| \, \max_p \sup_{[\alpha,\beta]}|\db{V_p}{x}|\,.
\end{equation}

Call $\xi_1 < \xi_2 < \cdots < \xi_{n-1}$ the zeros of the function $V_n$ in the interval $]\alpha,\beta[$.\medskip

\noib Assume that $V_n(\alpha) \not = 0$ and $V_n(\beta) \not = 0\,$.

Since $\db{V_n}{x}(\xi_i) \not = 0\,$, there exist $\delta_1, \varepsilon_1 > 0$ such that $|\db{V_n}{x}(x)| \ge \varepsilon_1$ for $ x \in [\xi_i - \delta_1,\xi_i+\delta_1]\,$, and $|V_n(x)| \ge \varepsilon_1$ in $[\alpha,\beta] \setminus \cup\, ]\xi_i-\delta_1,\xi_i+\delta_1[\,$.\\
 For $k$ large enough, we have $\omega M, \omega N \le \varepsilon_1 /2$. It follows that in the interval $[\xi_i-\delta_1,\xi_i+\delta_1]\,$,
$$ \left| \da (V_n + \omega \, \Pi)\right| \ge |\db{V_n}{x}| - \omega \, N \ge \varepsilon_1 /2\,.$$
 Furthermore,
$$
V_n(\xi_i \pm \delta_1) + \omega \, \Pi(\xi_i \pm \delta_1) \ge |V_n(\xi_i \pm \delta_1)| - \omega \, M \ge \varepsilon_1 /2\,.
$$
Since $V_n(\xi_i + \delta_1) V_n(\xi_i - \delta_1) < 0\,$, we can conclude that the function $V_n + \omega \, \Pi$ has exactly one zero in each interval  $]\xi_i-\delta_1,\xi_i+\delta_1[\,$.

In $[\alpha,\beta] \setminus \cup \, ]\xi_i-\delta_1,\xi_i+\delta_1[\,$, we have
$$
|V_n(x) + \omega \, \Pi(x)| \ge |V_n(x)| - \omega M \ge \varepsilon_1/2\,,
$$

which implies that $V_n(x) + \omega \, \Pi(x) \not = 0\,$.\medskip

\noib Assume that $V_n(\alpha) =0$ and $V_n(\beta) \neq 0\,$. This corresponds to the case $h=+\infty$ and $H \neq + \infty$. Hence the $V_j$ verify Dirichlet at $\alpha$ and $\Pi$ verifies Dirichlet at $\alpha$. Observing that $V'_n(\alpha) \neq 0$, it is immediate to see that there exists  $\delta_1>0$, such that, for $k$ large enough, $V_n(x) + \omega \, \Pi(x)$ has only $\alpha$ as zero in $[\alpha,\alpha +\delta_1]$. \medskip

\noib The other cases are treated in the same way. \hfill \qed

\clearpage

\section{Citations from Sturm's papers\\ French original and English translation}\label{S-Trans}

\textbf{Citation from \cite[Introduction]{Sturm1836a}}.\medskip

\begin{minipage}[t]{.45\linewidth}
On ne sait [ces \'{e}quations] les int\'{e}grer que dans un tr\`{e}s petit nombre de cas particuliers hors desquels on ne peut pas m\^{e}me en obtenir une int\'{e}grale premi\`{e}re ; et lors m\^{e}me qu'on poss\`{e}de l'expression de la fonction qui v\'{e}rifie une telle \'{e}quation, soit sous forme finie, soit en s\'{e}rie, soit en int\'{e}grales d\'{e}finies ou ind\'{e}finies, il est le plus souvent difficile de reconna\^{\i}tre dans cette expression la marche et les propri\'{e}t\'{e}s caract\'{e}ristiques de cette fonction. \ldots \\
S'il importe de pouvoir d\'{e}terminer la valeur de la fonction inconnue pour une valeur isol\'{e}e quelconque de la variable dont elle d\'{e}pend, il n'est pas moins n\'{e}cessaire de discuter la marche de cette fonction, ou en d'autres termes, d'examiner la forme et les sinuosit\'{e}s de la courbe dont cette fonction serait l'ordonn\'{e}e variable, en prenant pour abscisse la variable ind\'{e}pendante. Or on peut arriver \`{a} ce but par la seule consid\'{e}ration des \'{e}quations diff\'{e}rentielles elles-m\^{e}mes, sans qu'on ait besoin de leur int\'{e}gration.  Tel est l'objet du pr\'{e}sent m\'{e}moire. ~\ldots
\end{minipage}
\begin{minipage}[t]{0.5cm}\vrule \end{minipage}
\begin{minipage}[t]{.45\linewidth}
One only knows how to integrate these equations in a very small number of particular cases, and one can otherwise not even obtain a first integral; even when one knows the expression of the function which satisfies such an equation, in finite form, as a series, as integrals either definite or indefinite, it is most generally difficult to recognize in this expression the behaviour and the characteristic properties of this function. \ldots \\
Although it is important to be able to determine the value of the unknown function for an isolated value of the variable it depends upon, it is not less necessary to discuss the behaviour of this function, or otherwise stated, the form and the twists and turns of the curve whose ordinate would be the function, and the abscissa the independent variable. It turns out that one can achieve this goal by the sole consideration of the differential equation themselves, without having to integrate them. This is the purpose of the present memoir. ~\ldots
\end{minipage}

\newpage

\textbf{Citation from \cite[p.~108]{Sturm1836a}}.\medskip

\begin{minipage}[t]{.45\linewidth}
L'in\-t\'{e}\-grale compl\`{e}te de l'\'{e}quation (I) doit contenir deux constantes arbitraires, pour lesquelles on peut prendre les valeurs de $V$ et de $\db{V}{x}$ correspondantes \`{a} une valeur particuli\`{e}re de $x$. Lorsque ces valeurs sont fix\'{e}es, la fonction $V$ est enti\`{e}rement d\'{e}finie par l'\'{e}quation (I), elle a une valeur d\'{e}termin\'{e}e et unique pour chaque valeur de $x$.
\end{minipage}
\begin{minipage}[t]{1cm}\vrule \end{minipage}
\begin{minipage}[t]{.45\linewidth}
The complete integral of equation (I) must contain two arbitrary constants, for which one can take the values of $V$ and of $\db{V}{x}$ corresponding to some particular value of $x$. Once these values are fixed, the function $V$ is fully determined by equation (I), it has a uniquely determined value for each value of $x$.
\end{minipage}

\vspace{1cm}

\textbf{Citation from \cite[p.~379]{Sturm1836b}}.\medskip

\begin{minipage}[t]{.45\linewidth}
M. Liouville a d\'{e}montr\'{e} directement ce th\'{e}or\`{e}me, qui n'\'{e}tait pour moi qu'un corollaire du pr\'{e}c\'{e}dent, sans s'occuper du cas particulier o\`{u} la fonction serait nulle \`{a} l'une des extr\'{e}mit\'{e}s de la barre. J'en ai aussi trouv\'{e} apr\`{e}s lui une autre d\'{e}monstration directe que je donne dans ce m\'{e}moire. M.~Liouville a fait usage du m\^{e}me th\'{e}or\`{e}me dans un tr\`{e}s beau M\'{e}moire qu'il a publi\'{e} dans le num\'{e}ro de juillet de son journal et qui a pour objet le d\'{e}veloppement d'une fonction arbitraire en une s\'{e}rie compos\'{e}e de fonctions $V$ que nous avons consid\'{e}r\'{e}es.
\end{minipage}
\begin{minipage}[t]{1cm}\vrule \end{minipage}
\begin{minipage}[t]{.45\linewidth}
M. Liouville gave a direct proof of this theorem, which for me was a mere corollary of the preceding one, without taking care of the particular case in which the function vanishes at one of the extremities of the bar. I have also found, after him, another direct proof which I give in this memoir. M.~Liouville made use of the same theorem in a very nice memoir which he published in the July issue of his journal, and which deals with the expansion of an arbitrary function into a series made of the functions $V$ which we have considered.
\end{minipage}


\begin{weaker}
\clearpage
\begin{blue}
\section{Sturm's results under weaker assumptions}\label{S-weak}

We proved Theorems~\ref{T-st2} and \ref{T-LiSt} under the Assumptions~\eqref{E-intro-ass-S}. In this section, we consider the weaker assumptions
\begin{equation}\label{E-intro-ass-W}
\left\{
\begin{array}{l}
[\alpha,\beta] \subset ]\alpha_0,\beta_0[\,,\\[5pt]
K \in C^{1}(]\alpha_0,\beta_0[)\,,\\[5pt]
G, L \in C^{0}(]\alpha_0,\beta_0[)\,,\\[5pt]
K, G, L >0 \text{~on~} ]\alpha_0,\beta_0[\,.
\end{array}%
\right.
\end{equation}

Under these assumptions, the functions $V_j$ are $C^{2}$ on $]\alpha_0,\beta_0[$.  This follows easily for example from Liouville's existence proof \cite{Liou1836a}, and we have the following lemma, whose proof is analogous to the proof of Lemma~\ref{L-stp-2}

\begin{lemma}\label{L-stpw-2}
Let $k \in \Z\,$.
\begin{enumerate}
  \item The function $Y_k$ satisfies the boundary conditions \eqref{E-bca} and \eqref{E-bcb}.
  \item The functions $Y_k$ and $Y_{k+1}$ satisfy the relation
\begin{equation}\label{E-stp-6k0}
G\, Y_{k+1} = K\, \ddb{Y_k}{x} + \db{K}{x}\, \db{Y_k}{x} - L\, Y_k\,.
\end{equation}
  \item Under the Assumptions~\eqref{E-intro-ass-W}, the function $Y_k$ cannot vanish identically on an open interval $]\alpha_1,\beta_1[ \subset ]\alpha_0,\beta_0[$, unless $Y \equiv 0\,$.
\end{enumerate}
\end{lemma}%

In Subsection~\ref{SS-lip-1}, we have used Lemma~\ref{L-lip-2}~(2) which relies on the fact that the functions $V_j$ are $C^{\infty}$. If the functions $V_j$ are only $C^2$, we can apply Lemma~\ref{L-lip-2}~(1). It is easy to conclude that Liouville's proofs in Subsection~\ref{SS-lip-1} and \ref{SS-lip-2} go through, under the weaker Assumptions~\eqref{E-intro-ass-W}, if we only count distinct zeros, see \eqref{E-stp-16}. More precisely, we can prove the following claim.

\begin{claim}\label{C-lipT2}
Under the Assertions~\eqref{E-intro-ass-W}, for any $k \in \Z\,$, if the function $Y_k$ has at least $\mu$ distinct zeros in the interval $]\alpha,\beta[\,$, then the function $Y_{k+1}$ has at least $\mu$ distinct zeros in the interval $]\alpha,\beta[\,$.
\end{claim}%

We can then deduce from this claim, as in Section~\ref{S-lip}, that a linear combination $Y = \sum_{j=m}^n A_j V_j$ has at most $(n-1)$ distinct zeros (in particular it has finitely many zeros).

Once this result is secured, we can define zeros at which $Y$ changes sign (without using the multiplicity), and apply Sturm's lower bound argument to conclude that the function $Y$ must change sign at least $(m-1)$ times.
\medskip

\end{blue}%
\end{weaker}%

\begin{addendum}
\clearpage
\begin{blue}
\section{Sturm's original o.d.e proof}\label{S-sop}

The first proof of Theorem~\ref{T-st2r} appears in \cite[$\S$~{XXV}, p.~431]{Sturm1836b}, as a corollary of a more profound theorem ($\S$~{XXIV}) which describes the behaviour, as $t$ grows from $0$ to infinity, of the zeros of $x \mapsto u(x,t)$, where $u$ is a solution of the heat  \eqref{E-hist-2}-\eqref{E-hist-2ic}. \medskip

Sturm proves that the number $N(t)$ of zeros of the function $x \mapsto u(x,t)$ is piecewise constant, non-increasing in $t$, and that jumps occur precisely for values of $t$ such that $u(x,t)$ and $\frac{\partial u}{\partial t}(x,t)$ have common zeros. We refer to \cite{GaHa2005} for an analysis of this aspect of Sturm's paper \cite{Sturm1836b}.\medskip

The second proof, purely o.d.e., is  developed in \cite[$\S$~{XXVI}, p.~436 ff]{Sturm1836b}. In this section, we give the main steps of this proof (with page numbers and number of line from top $\ltop{}$, resp. from bottom $\lbot{}$).\medskip

\page{436} $\lbot{13}$, Sturm mentions Liouville's proof \cite{Liou1836b}.\\
{\it M. Liouville a d\'{e}montr\'{e} directement le th\'{e}or\`{e}me du num\'{e}ro pr\'{e}c\'{e}dent (dans le cahier d'ao\^{u}t de son journal) sans employer la consid\'{e}ration de la variable auxiliaire $t$ qui entre dans la fonction $u$ (42) dont j'ai fait usage. Il n'a pas tenu compte toutefois de la racine $\mathrm{x}$ ou $\mathrm{X}$\footnote{Respectively $\alpha$ and $\beta$ with our notation.} lorsqu'elle existe. Je vais donner ici une autre d\'{e}monstration directe du m\^{e}me th\'{e}or\`{e}me, ind\'{e}pendante de celui du n$^{\circ}${XXIV}.}\\
He introduces the linear combination
\begin{equation*}\label{E-43}
Y = C_i V_i + C_{i+1} V_{i+1} + \cdots + C_p V_p\,.\hspace{1cm} (43)
\end{equation*}
and, \page{436} $\lbot{1}$, its companion
\begin{equation*}\label{E-43a}
Y_1 = - \left( C_i \rho_i V_i + C_{i+1} \rho_{i+1} V_{i+1} + \cdots + C_p \rho_p V_p \right)\,.
\end{equation*}

\page{437}, Sturm establishes the differential relation
\begin{equation*}\label{E-44}
g Y_1 = k \ddb{Y}{x} + \db{k}{x} \db{Y}{x} - \ell Y\,. \hspace{1cm} (44)
\end{equation*}
He also notes $\ltop{5}$, that the function $Y$ satisfies the boundary conditions \eqref{E-bca}-\eqref{E-bcb}. Sturm' idea \emph{Je vais prouver \ldots }, is to prove that the function $Y_1$ has at least as many zeros in $]\alpha,\beta[$, counted with multiplicities, as the function $Y$ in the same circumstances.\medskip

\page{437} $\lbot{10}$, Sturm makes the \emph{implicit assumption} that the zeros of $Y$ are isolated. \medskip

\page{439} $\ltop{5}$, Sturm states that the number of sign changes of $Y_1$ in $]\alpha,\beta[$ is not smaller than the number of sign changes of $Y$. He then considers the zeros with multiplicities, and \emph{implicitly assumes} that the function $Y$ (assumed not to be identically zero) does not vanish at infinite order at some point. \medskip

\page{440} $\lbot{13}$, Sturm states that the number of zeros of $Y_1$ in $]\alpha,\beta[$, counted with multiplicities, is not smaller than the number of zeros of $Y$. He then examines ($\lbot{6}$) the possible zeros of $Y$ at $\alpha$ or $\beta$.\medskip

\page{442} $\ltop{7}$, Sturm states that the number of zeros of $Y_1$ in $[\alpha,\beta]$, counted with multiplicities (with a special rule for counting multiplicities at $\alpha$, $\beta$), is not smaller than the number of zeros of $Y$.

\page{442} $\ltop{13}$, Sturm iterates the procedure (with $Y_k$), and uses a limiting argument to conclude that the number of zeros of $Y$ in $[\alpha,\beta]$, counting multiplicities, is at most $p-1$. \medskip

\page{443}, Sturm proves the lower bound for the number of zeros and, ($\lbot{6}$), compares the present proof with the heat equation proof, the functions $Y_k$ are equal to $\frac{d^ku}{dt^k}(x,0)$. Finally, in a footnote, he mentions that $Y$ cannot vanish identically unless all the coefficients $C_j$ are zero. He does not mention the fact that $Y$ can actually not vanish at infinite order at any point.\medskip

\page{444} $\ltop{3}$, Sturm explains what to do when no assumption is made on the sign of the function $\ell$. Taking $Y$ as above, and defining
\begin{equation*}\label{E-43b}
Y_1 = - \left( C_i (\rho_i +c) V_i + C_{i+1} (\rho_{i+1} + c) V_{i+1} + \cdots + C_p (\rho_p + c) V_p \right)\,,
\end{equation*}
where $c$ is a constant, he obtains
\begin{equation*}\label{E-44a}
g Y_1 = k \ddb{Y}{x} + \db{k}{x} \db{Y}{x} - (g c + \ell) Y\,.
\end{equation*}
It suffices to assume that the constant $c$ is such that $gc + \ell > 0$ and to follow the previous proof with this new definition of $Y_1$.
\medskip


\section{Cross references to Sturm's and Liouville's papers}\label{S-cross}

In this Appendix, we give the references to pages in Sturm's paper \cite[$\S$~{XXVI}]{Sturm1836b} for the results in our paper.\medskip

\begin{itemize}
  \item Lemma~\ref{L-stp-2}: \page{437}. Note that the third assertion does not appear in Sturm's paper. He indeed implicitly assumes that the zeros of $Y$ are isolated.
  \item Lemma~\ref{L-stp-3}: \page{439}.
  \item Lemma~\ref{L-stp-4}: \page{440-441}.
  \item Lemma~\ref{L-stp-6}: \page{437}.
  \item Lemma~\ref{L-stp-8}: \page{438}.
  \item Proposition~\ref{P-stp-10}: \page{437-439}.
  \item Proposition~\ref{P-stp-12}: \page{439-442}.
  \item Proposition~\ref{P-stp-14}: \page{440-442}.
  \item Theorem~\ref{T-st2}: \page{442}.
\end{itemize}

Here are the pages in Liouville's paper \cite{Liou1836b}.

\begin{itemize}
  \item Theorem~\ref{T-LiSt}: \page{272}.
  \item Lemma~\ref{L-lip-2}: Mentioned \page{272}. No precise statement, no proof provided by Liouville.
  \item Proof of first assertion. Lemma~\ref{L-lip-4}: \page{274}.
  \item Proof of second assertion: \page{276} and reference to \cite{Liou1836a}.\\
      Claim~\ref{C-lip-2}: We use the determinant $\Delta$ to simplify Liouville's \cite[Lemme 1$^{\text{er}}$, p.~259]{Liou1836a}.
\end{itemize}
\medskip

\end{blue}%
\clearpage
\end{addendum}%

\newpage
\bibliographystyle{plain}
Numbers inserted after a reference indicate the pages where it is cited.

\end{document}